\documentclass{amsart}
\usepackage{amssymb}
\usepackage{amscd}
\usepackage[latin1]{inputenc}
\usepackage{amsmath}
\usepackage{amsfonts}
\usepackage{latexsym}
\usepackage{verbatim}
\usepackage{graphicx}
\usepackage{epsfig}
\newtheorem{theorem}{Theorem}
\newtheorem{corollary}[theorem]{Corollary}

\newtheorem{conj}[theorem]{Conjecture}

\theoremstyle{definition}

\newtheorem{remark}[theorem]{Remark}

\begin{document}
\title[Minimal zero entropy subshifts can be unrestricted along sparse sets]{Minimal zero entropy subshifts can be unrestricted along any sparse set}

\author{Ronnie Pavlov}
\address{Ronnie Pavlov\\
Department of Mathematics\\
University of Denver\\
2390 S. York St.\\
Denver, CO 80208}
\email{rpavlov@du.edu}
\urladdr{http://www.math.du.edu/$\sim$rpavlov/}

\thanks{The author gratefully acknowledges the support of a Simons Foundation Collaboration Grant.}
\renewcommand{\subjclassname}{MSC 2020}

\begin{abstract}
We present a streamlined proof of a result essentially present in \cite{PavlovThesis}, namely that 
for every set $S = \{s_1, s_2, \ldots\} \subset \mathbb{N}$ of zero Banach density and finite set $A$, there exists a minimal
zero-entropy subshift $(X, \sigma)$ so that for every sequence $u \in A^\mathbb{Z}$, there is $x_u \in X$ with $x_u(s_n) = u(n)$ for all $n \in \mathbb{N}$. Informally, minimal deterministic sequences can achieve completely arbitrary behavior upon restriction to a set of zero Banach density. 

As a corollary, this provides counterexamples to the Polynomial Sarnak Conjecture of \cite{polysarnak} which are significantly more general than some recently provided in \cite{KLR} and \cite{LianShi} and shows that no similar result can hold under only the assumptions of minimality and zero entropy.
\end{abstract}

\maketitle

\section{Introduction}\label{intro}

The well-known \textbf{Sarnak conjecture} states that the M\"{o}bius function $\mu$ is uncorrelated with all deterministic sequences. A sequence is called \textbf{deterministic} if it is the image under a continuous function of the trajectory of a point in a \textbf{topological dynamical system} with zero \textbf{entropy} (see Section~\ref{defs} for definitions of this and other concepts not defined in this introduction). More formally, 

\begin{conj}[Sarnak Conjecture]
If $(X,T)$ is a topological dynamical system with zero entropy, $x_0 \in X$, and $f \in C(X)$, then
\[
\frac{1}{N} \sum_{n=1}^N \mu(n) f(T^n x_0) \rightarrow 0.
\]
\end{conj}

Although this problem is still open, there are many recent works on the topic, which have made significant progress and resolved it for some classes of dynamical systems. 
In \cite{polysarnak}, a potential stronger `polynomial' (meaning that only polynomial iterates of $x_0$ are taken rather than all) version of the Sarnak Conjecture was conjectured. In order to rule out some degenerate examples, the assumption of \textbf{minimality} was added on $(X, T)$, meaning that for every $x \in X$, the set $\{T^n x\}$ is dense. 

\begin{conj}[Polynomial Sarnak Conjecture (\cite{polysarnak}, Conjecture 2.3)]\label{poly}
If $(X,T)$ is a minimal topological dynamical system with zero entropy, $x_0 \in X$, $f \in C(X)$, and $p: \mathbb{N} \rightarrow \mathbb{N}_0$ is a polynomial, then 
\[
\frac{1}{N} \sum_{n=1}^N \mu(n) f(T^{p(n)} x_0) \rightarrow 0.
\]
\end{conj}

This conjecture is now known to be false; recently Kanigowski, Lema\'{n}czyk, and Radziwi\l\l (\cite{KLR}) and Lian and Shi (\cite{LianShi}) have separately provided counterexamples. However, these counterexamples are specific to the case $p(n) = n^2$ (though they could perhaps be generalized) and make strong usage of the nice arithmetic properties of this function. The first is a skew product and the second is a symbolically defined dynamical system called a Toeplitz subshift. 

The purpose of this note is to show that even much weaker versions of Conjecture~\ref{poly} are false, because minimal zero entropy systems can achieve \textbf{any} possible behavior (i.e., not just correlation with $\mu$) along \textbf{any} prescribed set $S \subset \mathbb{N}$ of zero \textbf{Banach density} (i.e., not just the image of a polynomial). One such result had already been proved by the author in \cite{PavlovThesis}, which already immediately refutes the Polynomial Sarnak Conjecture. 

\begin{theorem}[\cite{PavlovThesis}, Corollary 5.1]\label{oldthm}
Assume that $d \in \mathbb{N}$, $(w_n)$ is an increasing sequence of positive integers where $w_{n+1} < (w_{n+1} - w_n)^{d+1}$ for large enough $n$, and $(z_n)$ is any sequence in $\mathbb{T} := \mathbb{Z}/\mathbb{N}$. Then there exists a totally minimal, totally uniquely ergodic, topologically mixing zero entropy map $S$ on $\mathbb{T}^{2d+4}$ so that, if $\pi$ is projection onto the final coordinate, $\pi(S^{w_n} \mathbf{0}) = z_n$ for sufficiently large $n$.
\end{theorem}

(We don't further work with the properties of unique ergodicity and topological mixing, and so don't provide definitions here. However, we do note that Theorem~\ref{oldthm} shows that even adding these hypotheses to Conjecture~\ref{poly} would not make it true.)
We note that the entropy of the transformation $S$ was never mentioned in \cite{PavlovThesis}. However, $S$ is defined as a suspension flow of a product of a toral rotation and a skew product $T$ under a roof function $1 < g < 3$. The skew product $T$ is of the form $(x_1, x_2, x_3, \ldots, x_m) \mapsto (x_1 + \alpha, x_2 + f(x_1), x_3 + x_2, \ldots, x_m + x_{m-1})$ for a continuous self-map $f$ of $\mathbb{T}$. Since its first coordinate is an irrational rotation, known to have zero entropy, the map $T$ also has zero entropy by Abramov's skew product entropy formula. Then $S$ has zero entropy as well, by Abramov's suspension flow entropy formula.

\begin{remark} 
Here are a few more relevant facts about the construction from \cite{PavlovThesis}:
\begin{enumerate}
\item The map $S$ is \textbf{distal}, meaning that for all $x \neq y$, $\{d(T^n x, T^n y)\}_n$ is bounded away from $0$.
\item The roof function $g$ is $C^{\infty}$ and the function $f$, though not $C^\infty$ as constructed in \cite{PavlovThesis}, can easily be made so; it is just a uniformly convergent infinite series of `bump functions,' which can easily be chosen $C^\infty$. 
\end{enumerate}
The second fact may be of interest since the authors of \cite{KLR} prove a positive result for convergence along prime iterates of similar skew products $(x, y) \mapsto (x+\alpha, y+f(x))$ under the assumption that the function $f$ is real analytic, provide some counterexamples with continuous $f$, and ask whether this assumption could be weakened to $C^\infty$. Though the constructions are not exactly the same, and though the primes absolutely do not satisfy the assumption of Theorem~\ref{oldthm}, (2) might suggest that $C^\infty$ is not always sufficient for good averaging of skew products along sparse sequences.
\end{remark}

We note that Theorem~\ref{oldthm} clearly applies to any sequence $w_n = p(n)$ for a nonconstant polynomial $p: \mathbb{N} \rightarrow \mathbb{N}_0$ (possibly omitting finitely many terms), and so, by simply defining $z_n$ to be $\frac{1}{2}$ when $\mu(n) = 1$ and $0$ otherwise, one achieves
\[
\frac{1}{N} \sum_{n=1}^N \mu(n) \pi(S^{p(n)} \textbf{0}) = 
\frac{0.5 |\mu^{-1}(\{1\}) \cap \{1, \ldots, N\}|}{N},
\]
which does not approach $0$ as $N \rightarrow \infty$, disproving the Polynomial Sarnak Conjecture for every nonconstant $p$. The same is true of any function $p$ with polynomial growth, even for degree less than $2$, e.g. $p(n) = \lfloor n^{1.01} \rfloor$. However, Theorem~\ref{oldthm} does not apply to more slowly growing $p$ such as $\lfloor n \ln n \rfloor$. The author proved a different result (Corollary 3.1) in \cite{PavlovThesis} using \textbf{subshifts}; a subshift is a closed shift-invariant subset of 
$A^{\mathbb{Z}}$ (for some finite alphabet $A$) endowed with the left-shift transformation. 
Corollary 3.1 of \cite{PavlovThesis} states that given any sequence of zero Banach density (regardless of growth rate), there exists a minimal subshift whose points can achieve arbitrary behavior along that sequence. However, entropy was not mentioned there, and although the proof there can indeed yield a zero entropy subshift, it's not easy to verify; the construction is quite complicated in order to achieve $(X,T)$ which is totally minimal, totally uniquely ergodic, and topologically mixing.

In this note, we present a streamlined self-contained proof of the following result, which shows that minimal zero entropy subshifts can realize arbitrary behavior along any sequence of zero Banach density.

\begin{theorem}\label{mainthm}
For any $S = \{s_1, s_2, \ldots\} \subset \mathbb{N}$ with $d^*(S) = 0$ and any finite alphabet $A$, there exists a minimal zero entropy subshift $X \subset A^\mathbb{Z}$ so that for every $u \in A^\mathbb{N}$, there is 
$x_u \in X$ where $x_u(s_n) = u(n)$ for all $s \in S$.
\end{theorem}

We note that this proves that even with substantially weaker hypotheses, nothing in the spirit of the Polynomial Sarnak Conjecture can hold under only the assumptions of minimality and zero entropy. Even if $p$ is only assumed to have range of zero Banach density and $\rho: \mathbb{N} \rightarrow \mathbb{Z}$ is only assumed to have $\limsup \frac{1}{N} \sum_{n=1}^N |\rho(n)| > 0$ (equivalently, $\rho$ takes nonzero values on a set of positive upper density), one can define a subshift $X$ on $\{-1,0,1\}$ and $x_u \in X$ as in Theorem~\ref{mainthm} for $u(n) = \textrm{sgn}(\rho(n))$. 
Then, for $f \in C(X)$ defined by $x \mapsto x(0)$, 
the limit supremum of the averages
\begin{multline*}
\frac{1}{N} \sum_{n=1}^N \rho(n) f(\sigma^{p(n)} x_u) = 
\frac{1}{N} \sum_{n=1}^N \rho(n) x_u(p(n)) = 
\frac{1}{N} \sum_{n=1}^N \rho(n) u(n) = \\
\frac{1}{N} \sum_{n=1}^N \rho(n) \textrm{sgn}(\rho(n)) = 
\frac{1}{N} \sum_{n=1}^N |\rho(n)|
\end{multline*}
is positive by assumption.

We remark that when $\rho = \mu$ is the M\"{o}bius function, this means that 
\[
\frac{1}{N} \sum_{n=1}^N \mu(n) f(\sigma^{p(n)} x_u)
\]
can be made to approach $\frac{6}{\pi^2}$ (for $x_u$ in a minimal zero-entropy subshift), a slight improvement of \cite{LianShi} which showed that it could attain values arbitrarily close to $\frac{6}{\pi^2}$. 

\section{Definitions}\label{defs}

A \textbf{topological dynamical system} $(X, T)$ is defined by a compact metric space $X$ and homeomorphism $T: X \rightarrow X$. A \textbf{subshift} is a topological dynamical system defined by some finite set $A$ (called the \textbf{alphabet}) and the restriction of the \textbf{left shift} map $\sigma: A^{\mathbb{Z}} \rightarrow A^{\mathbb{Z}}$ defined by $(\sigma x)(n) = x(n+1)$ to some closed and $\sigma$-invariant $X \subset A^{\mathbb{Z}}$ (with the induced product topology). A subshift $(X, \sigma)$ is \textbf{minimal} if for every $x \in X$, $\{\sigma^n x\}_{n \in \mathbb{Z}}$ is dense in $X$.

A \textbf{word} over $A$ is any finite string of symbols from $A$; a word $w = w(1) \ldots w(n)$ is said to be a \textbf{subword} of a word or infinite sequence $x$ if there exists $i$ so that $w(1) \ldots w(n) = x(i+1) \ldots x(i+n)$. 
The \textbf{language} $L(X)$ of a subshift $(X, \sigma)$ is the set of all subwords of sequences in $X$, and for any $n \in \mathbb{N}$ we denote $L_n(X) = L(X) \cap A^n$. For two words $u = u(1) \ldots u(m)$ and $v = v(1) \ldots v(n)$, denote by $uv$ their \textbf{concatenation} $u(1) \ldots u(m) v(1) \ldots v(n)$.

We do not give a full definition of \textbf{topological entropy} here, but note that it is a number $h(X,T) \in [0, \infty]$ associated to any TDS $(X, T)$ which is conjugacy-invariant. We will only need the following definition for subshifts: for any $(X, \sigma)$, 
\[
h(X, \sigma) = \lim \frac{\ln |L_n(X)|}{n}.
\]

The \textbf{Banach density} of a set $S \subset \mathbb{N}$ is
\[
d^*(S) := \lim_{n \rightarrow \infty} \sup_{k \in \mathbb{N}} \frac{|S \cap \{k, \ldots, k+n-1\}|}{n}.
\]

\section{Proof of Theorem~\ref{mainthm}}\label{proof}

\begin{proof}

As in \cite{PavlovThesis}, we adapt the block-concatenation construction of Hahn and Katznelson (\cite{HK}).

We construct $X$ iteratively via auxiliary sequences $m_k$ of odd positive integers, $A_k \subset A^{m_k}$, and $w_k \in A_k$. 
Define $m_0 = 1$, $A_0 = A$, and $w_0 = 0$ (which we assume without loss of generality to be in $A$). Now, suppose that $m_k$, $A_k$, and $w_k$ are defined. Define 
$m_{k+1} > \max(3m_k |A_k|, 12(\ln 2)(4/3)^{k+1})$ to be an odd multiple of $3m_k$ large enough that $|S \cap I|/|I| < (3m_k)^{-1}$ for all intervals $I$ of length $m_{k+1}$ (using the fact that $d^*(S) = 0$). Define $A_{k+1}$ to be the set of all concatenations of $\frac{m_{k+1}}{m_k}$ words in $A_k$
in which every word in $A_k$ is used at least once and in which at least one-third of the concatenated words are equal to $w_k$. 
Define $Y_k$ to be the set of shifts of biinfinite (unrestricted) concatenations of words in $A_k$, define $Y = \bigcap_k Y_k$, and define
$X$ to be the subshift of $Y$ consisting of sequences in which every subword is a subword of some $w_k$.

We claim that $(X, \sigma)$ is minimal. Indeed, consider any $x \in X$ and $w \in L(X)$. By definition, $w$ is a subword of $w_k$ for some $k$. By definition, $w_k$ is a subword of every word in $A_{k+1}$. Finally, $x$ is a shift of a concatenation of words in 
$A_{k+1}$, each of which contains $w_k$, and therefore $w$. So, $x$ contains $w$, and since $w \in L(X)$ was arbitrary, the orbit of $x$ is dense. Since $x \in X$ was arbitrary, \textbf{$(X, \sigma)$ is minimal}. 

We also claim that $(X, \sigma)$ has zero entropy. We see this by bounding $|A_k|$ from above. For every $k$, each word in $A_{k+1}$ is defined by
an ordered $(m_{k+1}/m_k)$-tuple of words in $A_k$, where at least one-third are $w_k$. The number of such tuples can be bounded from above by
\[
{m_{k+1}/m_k \choose m_{k+1}/3m_k} |A_k|^{2m_{k+1}/3m_k} \leq 2^{m_{k+1}/m_k} |A_k|^{2m_{k+1}/3m_k}.
\]
Therefore,
\[
\frac{\ln |A_{k+1}|}{m_{k+1}} \leq \frac{\ln 2}{m_k} + \frac{2}{3} \frac{\ln |A_k|}{m_k}.
\]
Now, it's easily checked that $\frac{\ln |A_k|}{m_k} \leq \ln |A| (3/4)^k$ for all $k$ by induction. The base case $k = 0$ is immediate. For the inductive step, if we assume that 
$\frac{\ln |A_k|}{m_k} \leq \ln |A| (3/4)^k$, then recalling that $m_k > 12(\ln 2) (4/3)^k$,
\[
\frac{\ln |A_{k+1}|}{m_{k+1}} < \frac{1}{12} (3/4)^{k} + \frac{2}{3} \ln |A| (3/4)^k \leq \frac{\ln |A|}{12} (3/4)^{k} + \frac{2}{3} \ln |A| (3/4)^k = \ln |A| (3/4)^{k+1}.
\]
Therefore, for all $k$, $|A_k| \leq e^{\ln |A| (3/4)^k m_k}$. Finally, we note that every word in $L_{m_k}(X)$ is a subword of a concatenation of a pair of words in $A_k$, so determined by such a pair and by the location of the first letter. Therefore, $|L_{m_k}(X)| \leq m_k |A_k|^2 < m_k e^{2\ln |A|(3/4)^k m_k}$. This clearly implies that
\[
h(X) = \lim_{k \rightarrow \infty} \frac{\ln |L_{m_k}(X)|}{m_k} \leq \limsup_{k \rightarrow \infty} \frac{\ln m_k}{m_k} + 2\ln |A|(3/4)^k = 0,
\]
i.e. \textbf{$\mathbf{X}$ has zero entropy}.

It remains, for $u \in A^\mathbb{N}$, to construct $x_u \in X$ with $x_u(s_n) = u(n)$ for all $s_n \in S$. The construction of $x_u$ proceeds in steps, where it is continually assigned letters from $A$ on portions of $\mathbb{Z}$, with undefined portions labeled by $*$. Formally, define $x^{(0)} \in A \sqcup \{*\}^{\mathbb{Z}}$ by $x^{(0)}(s_n) = u(n)$ for $s \in S$ and $*$ for all other locations. 

Now partition $\mathbb{Z}$ into the intervals $((i-0.5)m_1, (i+0.5)m_1)$ (herein, all intervals are assumed to be intersected with $\mathbb{Z}$).
For every $i$ for which $S \cap ((i-0.5)m_1, (i+0.5)m_1) \neq \varnothing$, consider the $m_1$-letter word
$x^{(0)}(((i-0.5)m_1, (i+0.5)m_1))$. 
By definition of $m_1$, $|S \cap ((i-0.5)m_1, \ldots, (i+0.5)m_1)| < m_1/3m_0 = m_1/3$, and so at most one-third
of the letters in this word are non-$*$. Fill the remaining locations by assigning the first $m_1/3$ as
$w_0 = 0$. At least $m_1/3$ letters remain, which is larger than $|A_0| = |A|$ by definition of $m_1$.
Fill those in an arbitrary way which uses all letters from $A$ at least once. The resulting $m_1$-letter word
is in $A_1$ by definition, call it $w^{(1)}_i$. Now, define $x^{(1)}$ by setting
$x^{(1)}(((i-0.5)m_1, (i+0.5)m_1)) = w^{(1)}_i$ for all $i$ as above (i.e. those for which $S \cap ((i-0.5)m_1, (i+0.5)m_1) \neq \varnothing$)
and $*$ elsewhere. 
Note that $x^{(1)}$ is an infinite concatenation of words in $A_1$ and blocks of $*$ of length $m_1$ and that $x^{(1)}$ contains $*$ on any interval $((i-0.5)m_1, (i+0.5)m_1)$ which is disjoint from $S$.

Now, suppose that $x^{(k)}$ has been defined as an infinite concatenation of words in $A_k$ and blocks of $*$ of length $m_k$ which contains $*$ on any interval $((i-0.5)m_k, (i+0.5)m_k)$ which is disjoint from $S$. We wish to extend $x^{(k)}$ to $x^{(k+1)}$ by changing some $*$ symbols to letters in $A$. Consider any
$i$ for which $S \cap ((i-0.5)m_{k+1}, \ldots, (i+0.5)m_{k+1}) \neq \varnothing$. The portion of $x^{(k)}$ occupying that interval is a concatenation of words in $A_k$ and blocks of $*$ of length $m_k$ (we use here the fact that $m_{k+1}$ is odd), and the number which are words in $A_k$ is bounded from above by
the number of $j \in ((i-0.5)m_{k+1}/m_k, (i+0.5)m_{k+1}/m_k)$ for which $((j-0.5)m_k, (j+0.5)m_k)$ is not disjoint from $S$, which in turn is bounded from above by $|S \cap ((i-0.5)m_{k+1}, (i+0.5)m_{k+1})|$, which by definition of $m_{k+1}$ is less than $m_{k+1}/3m_k$. 
Therefore, at least two-thirds of the concatenated $m_k$-blocks comprising $x^{(k)}(((i-0.5)m_{k+1}, (i+0.5)m_{k+1}))$ are blocks of $*$. 
Fill the first $m_{k+1}/3m_k$ of these with $w_k$. Then at least $m_{k+1}/3m_k$ blocks remain, which is more than $|A_k|$ by definition of $m_{k+1}$. Fill these in an arbitrary way which uses each word in $|A_k|$ at least once. By definition, this creates a word in $A_{k+1}$, which we denote by
$w^{(k+1)}_i$. Define $x^{(k+1)}(((i-0.5)m_{k+1}, (i+0.5)m_{k+1})) = w^{(k+1)}_i$ for any $i$ as above (i.e. those for which $S \cap ((i-0.5)m_{k+1}, (i+0.5)m_{k+1}) \neq \varnothing$) and as $*$ elsewhere. Note that $x^{(k+1)}$ is an infinite concatenation of words in 
$A_{k+1}$ and blocks of $*$ of length $m_{k+1}$ which contains $*$ on any interval $((i-0.5)m_{k+1}, (i+0.5)m_{k+1})$ which is disjoint from $S$.

We now have defined $x^{(k)} \in (A \sqcup \{*\})^{\mathbb{Z}}$ for all $k \in \mathbb{N}$. Since each is obtained from the previous by changing some $*$s to letters from $A$, they approach a limit $x_u$ which agrees with $x^{(0)}$ on all locations where $x^{(0)}$ had 
letters from $A$, i.e. $x_u(s_n) = u(n)$ for all $n \in \mathbb{N}$. Since $S \neq \varnothing$, $S \cap (-0.5m_k, 0.5m_k) \neq \varnothing$ for all large enough $k$, and so $x^{(k)}((-0.5m_k, 0.5m_k))$ has no $*$, meaning that $x_u \in A^\mathbb{Z}$. 

It remains only to show that $x_u \in X$. By definition, $x_u$ is a concatenation of words in $A_k$ for every $k$, so $x_u \in Y = \bigcap_k Y_k$ as in the definition of $X$. Finally, every subword $w$ of $x_u$ is contained in $x_u((-0.5m_k, 0.5m_k))$ for large enough $k$, and this word is in $A_k$ by definition. Since all words in $A_k$ are subwords of $w_{k+1}$, $w$ is also. Therefore by 
definition, $\mathbf{x_u \in X}$ and $\mathbf{x_u(s_n) = u(n)}$ for all $n$, completing the proof.

\end{proof}

\begin{remark}

We observe that the assumption of zero Banach density cannot be weakened in Theorem~\ref{mainthm}. Assume for a contradiction that $S \subset \mathbb{N}$ has $d^*(S) = \alpha > 0$, and that every $u \in A^\mathbb{N}$ could be assigned $x_u$ as in Theorem~\ref{mainthm}. By definition of Banach density, there exist intervals $I_n$ with lengths approaching infinity so that $|S \cap I_n|/|I_n| > \alpha/2$ for all $n$. For every $n$, since all possible assignments of letters from $A$ to locations in $S \cap I_n$ give rise to sequences in $X$,
$|L_{|I_n|}(X)| \geq 2^{|S \cap I_n|} > |A|^{\alpha|I_n|/2}$. Then,
\[
h(X) = \lim_n \frac{\ln |L_{|I_n|}(X)|}{|I_n|} \geq \limsup \frac{\ln |A|^{\alpha|I_n|/2}}{|I_n|} = \alpha (\ln |A|)/2 > 0.
\]
Therefore, no such $X$, minimal or otherwise, can have zero entropy.
\end{remark}

\bibliographystyle{plain}
\bibliography{sarnak}

\end{document}